# A Jump Ornstein–Uhlenbeck Bridge Based on Energy-optimal Control and Its Self-exciting Extension

Hidekazu Yoshioka, *Member, IEEE* and Kazutoshi Yamazaki

*Abstract*— We study a version of the Ornstein–Uhlenbeck bridge driven by a spectrally-positive subordinator. Our formulation is based on a Linear-Quadratic control subject to a singular terminal condition. The Ornstein–Uhlenbeck bridge, we develop, is written as a limit of the obtained optimally controlled processes, and is shown to admit an explicit expression. Its extension with self-excitement is also considered. The terminal condition is confirmed to be satisfied by the obtained process both analytically and numerically. The methods are also applied to a streamflow regulation problem using a real-life dataset.

## I. Introduction

### A. Study Background

The Ornstein-Uhlenbeck (OU) process is a stochastic process with applications in various fields ranging from mathematical finance to physical sciences. It is a mean-reverting process which tends to drift toward a certain value (mean), particularly applicable in modelling asset prices in financial markets. The most classical and well-known OU process is the one driven by standard Brownian motion. However, because the OU process – which belongs to the class of affine processes – has various analytical properties, its generalizations and extensions are often tractable. The continuous-state branching process with immigration (CBI process) [1] can be seen as an extension of the OU process with self-exciting behaviors. These processes have been applied in a variety of industrial problems such as currency option pricing [2], evaluation of power markets [3], and water environmental assessment [4].

In this paper, we are interested in a (stochastic) bridge version of the OU process, or a modification of the OU process obtained by imposing a constrained terminal condition so that it ends at some given value at some prespecified time. The most classical stochastic bridge is the Brownian bridge, or a Brownian motion conditioned to satisfy the terminal constraint [5], and there also exist several examples of bridges obtained from other Gaussian processes. In general, a stochastic bridge is not unique and there are several ways to construct it. For example, Doob's $h$-transform is used to construct a Brownian bridge [6], whereas a Brownian-motion-driven OU bridge was obtained through a stochastic control approach [7]. The process obtained in the latter has been applied in a wide variety of research fields owing to their analytical tractability, such as statistical inference [8], commodity pricing [9], error analysis of stochastic partial differential equations [10], energy demand modeling [11], and quasi-Monte-Carlo methods [12]. A time-fractional version of the OU bridge has also been discussed in [13]. The OU bridge under additional constraints regarding the integrated value of the sample path has been analyzed by [14]. However, these papers on OU bridges and a majority in the existing literature focus on Gaussian cases; results on non-Gaussian OU bridges containing jumps are limited. Further, an extension driven by self-exciting processes, as we consider in this paper, has not been studied.

### B. Objective and Contribution

The objective of this paper is to propose and analyze a subordinator-driven OU bridge. Contrary to the classical OU bridge driven by a Brownian motion, our version admits jumps in the process. This paper focuses on the process driven by a pure-jump subordinator with non-Gaussian noise, relevant to the mathematical modelling of the streamflow discharge such as flood control.

To construct our OU bridge, we impose the terminal condition based on a certain finite-horizon linear-quadratic (LQ) regulator as in Chen and Georgiou [7]. More specifically, an LQ control problem minimizing the expected quadratic control cost and terminal cost ($L^2$ norm of the difference between the process at the terminal time and the target) is firstly formulated, where the penalty of the terminal cost is modulated by a parameter. The OU bridge is obtained by taking the parameter to infinity. Note that this energy-minimization-based construction of a stochastic bridge is not, in general, equivalent to that obtained by applying Doob's $h$-transform [14], which relies heavily on the closed-form representation of some conditional probabilities that are not available for a general process with jumps. By contrast, our approach via LQ control is potentially applicable in a wider setting driven by other (affine) processes.

Because our OU bridge is written as a limit of the optimally controlled processes, precise arguments are necessary to confirm that the obtained process indeed satisfies the terminal condition almost surely (a.s.). This is verified by analyzing the growth speed of the control process as it approaches the terminal time. Similar control problems subject to (possibly stochastic) terminal constraints have been studied from the standpoint of the dynamic programming and backward stochastic differential equations [15-16] with some applications in economics [17-18]. However, the bridges of the form developed in this paper have not been studied. Various quantities of interest, including moments, can be computed efficiently by solving ordinary differential equation (ODE) systems obtained in this paper. We further consider a

*Research supported by KAKENHI No. 22K14441 and 22H02456 from Japan Society for the Promotion Science and the start-up grant by the School of Mathematics and Physics of the University of Queensland. A part of this work was conducted during the first author was at Shimane University.

H. Yoshioka is with Japan Advanced Institute of Science and Technology, 1-1 Asahidai, Nomi, Ishikawa 923-1292, Japan (Tel: +81-761-51-1745; e-mail: yoshih@jaist.ac.jp). K. Yamazaki is School of Mathematics and Physics, The University of Queensland, St Lucia, Brisbane, QLD 4072, Australia (e-mail: k.yamazaki@uq.edu.au).

generalization of the OU bridge having a self-exciting jump term, which we call the self-exciting (SE) bridge.

Computational experiments based on a Monte-Carlo method are further conducted to demonstrate that our theoretical results agree with those computed by Monte Carlo simulation for both OU and SE bridges. As their applications in the streamflow discharge regulation problem, we calibrate the SE bridge using real-life discharge data and compute an optimal streamflow when the discharge needs to be maintained to be at a level at given terminal time.

## II. ORNSTEIN–UHLENBECK BRIDGE

### A. Ornstein–Uhlenbeck Process

On a probability space $(\Omega, \mathcal{F}, \mathbb{P})$, we define an OU process $X$ as a solution to the linear stochastic differential equation (SDE)

$$dX_t = -rX_t dt + \int_0^{+\infty} zN(dz, dt) \text{ for } t > 0, \ X_0 = x_0 \in \mathbb{R}. \quad (1)$$

Here, $r > 0$ is a reversion rate and $N$ is a Poisson random measure associated with a (spectrally-positive) subordinator whose Lévy measure $v$ on $(0,\infty)$ admits the first two moments, i.e., $M_k := \int_0^{+\infty} z^k v(dz) < +\infty$ ($k = 1, 2$). The SDE (1) has a unique strong càdlàg solution defined globally in time (Theorem 4.76 of [19]).

### B. Linear-Quadratic Control Problem

We construct a stochastic bridge in terms of (the limiting case of) an LQ control problem. We assume a fixed terminal time $T = 1$ without loss of generality; cases with generic $T > 0$ can be handled by properly scaling parameters. The objective is to control the paths of the OU process (1) to satisfy the terminal constraint. In other words, we want to construct a modification of $X$, say $X^*$, such that

$$X_1^* = \hat{x} \text{ a.s.} \quad (2)$$

with a prescribed deterministic constant $\hat{x} \in \mathbb{R}$. Following [7], we consider a parameterized LQ control problem and take their certain limit to enforce the terminal condition (2).

An admissible control $u$ is an $\mathbb{R}$-valued process adapted to the filtration generated by the OU process $X$ with $\mathbb{E}\left[\int_0^1 u_s^2 ds\right] < +\infty$. It modifies the drift term of (1) and the corresponding controlled process $X^{(u)}$ becomes a solution to the SDE

$$dX_t^{(u)} = \left(-rX_t^{(u)} + u_t\right)dt + \int_0^{+\infty} zN(dz, dt) \text{ for } t > 0, \quad (3)$$

starting at $X_0^{(u)} = x_0$. We consider the minimization problem

$$\inf_u \mathbb{E}\left[\frac{F}{2}\left(X_1^{(u)} - \hat{x}\right)^2 + \int_0^1 \frac{1}{2}u_s^2 ds\right] \quad (4)$$

over all admissible controls $u$. Here, $F > 0$ is a penalization parameter of the terminal condition.

We solve the control problem (4) for each $F > 0$ and then consider the limit of the optimally controlled process as $F \to +\infty$, for which the constraint (2) is confirmed to hold. Fix $F > 0$. The solution to the control problem (4) can be explicitly written in terms of the solution $\Phi(t, x)$ ($0 \le t \le 1$, $x \in \mathbb{R}$) to the following non-local Hamilton–Jacobi–Bellman (HJB) equation

$$\frac{\partial \Phi}{\partial t} - rx\frac{\partial \Phi}{\partial x} + \int_0^{+\infty}\left(\Phi(t, x+z) - \Phi(t, x)\right)v(dz) + \inf_{u \in \mathbb{R}}\left\{u\frac{\partial \Phi}{\partial x} + \frac{1}{2}u^2\right\} = 0 \quad (5)$$

for $0 < t < 1$ and $x \in \mathbb{R}$ subject to the terminal condition

$$\Phi(1, x) = \frac{F}{2}(x - \hat{x})^2 \text{ for } x \in \mathbb{R}. \quad (6)$$

After solving this, the solution to (4) is given by $\Phi(0, x_0)$.

To solve (5)-(6), we guess a solution of the form

$$\Phi(t, x) = \frac{1}{2}A_{t,F}x^2 + B_{t,F}x + C_{t,F} \text{ for } 0 \le t \le 1 \text{ and } x \in \mathbb{R} \quad (7)$$

for time-dependent deterministic coefficients $A_{t,F}, B_{t,F}, C_{t,F}$ parameterized by $F$. A direct substitution of (7) in (5)-(6) yields the Riccati equation given by, for $0 \le t < 1$,

$$\frac{dA_{t,F}}{dt} = A_{t,F}^2 + 2rA_{t,F}, \quad (8)$$

$$\frac{dB_{t,F}}{dt} = (r + A_{t,F})B_{t,F} - M_1 A_{t,F}, \quad (9)$$

$$\frac{dC_{t,F}}{dt} = \frac{1}{2}B_{t,F}^2 - M_1 B_{t,F} - \frac{1}{2}M_2 A_{t,F} \quad (10)$$

subject to the terminal conditions

$$A_{1,F} = F, \ B_{1,F} = -F\hat{x}, \ C_{1,F} = F\hat{x}^2. \quad (11)$$

Hereafter, we set $q := F(F + 2r)^{-1} \in (0, 1)$. The optimal control, which minimizes (4) when $X_t^{(u)} = x$ is given as

$$u_F^*(t, x) := \arg\min_u \left\{u\frac{\partial \Phi}{\partial x} + \frac{u^2}{2}\right\} = -\left(A_{t,F}x + B_{t,F}\right) \quad (12)$$

for $0 \le t \le 1$ and $x \in \mathbb{R}$, with the solutions $A_{t,F}, B_{t,F}$ to the Riccati equation given by

$$A_{t,F} = \frac{2rq}{\exp(2r(1-t)) - q}, \quad (13)$$

$$B_{t,F} = -F\hat{x}\frac{1-q}{\exp(r(1-t))(1 - q\exp(-2r(1-t)))} + \frac{2qM_1(1 - q\exp(-r(1-t)))}{\exp(r(1-t))(1 - q\exp(-2r(1-t)))}. \quad (14)$$

The optimality of the control given by (12) for the problem (4) can be verified rigorously following a classical argument since the solution (7) to the HJB equation (5) is smooth and grows only quadratically (see, e.g., [Theorem 5.1 of 20]).

We now take the penalty parameter $F \to +\infty$ and study the limit $u^*$ of the optimal control $u_F^*$. Letting $F \to +\infty$ (and hence $q \to 1$) in (13)-(14) yields the following convergence, uniform in $t$ on any compact subset of $[0,1)$:

$$\lim_{F \to +\infty} A_{t,F} = A_t := \frac{2r}{\exp(2r(1-t))-1}, \quad (15)$$

$$\lim_{F \to +\infty} B_{t,F} = B_t := \frac{2r\left(-\hat{x} + \frac{M_1}{r}\left(1-\exp(-r(1-t))\right)\right)}{\exp(r(1-t))\left(1-\exp(-2r(1-t))\right)}. \quad (16)$$

As $t \to 1-$, $A_t, B_t = O\left((1-t)^{-1}\right)$, which is a speed sufficient to satisfy the terminal condition, as we confirm in the next section. With $u_F^*$ defined by (12) with $A_{t,F}$ and $B_{t,F}$ replaced by (15) and (16), respectively, we consider the process $X_t^*$ that solves the SDE

$$dX_t^* = \left(-(r+A_t)X_t^* - B_t + M_1\right)dt + \int_0^{+\infty} z\tilde{N}(dz,dt), \quad (17)$$

where $\tilde{N}(dz,ds) = N(dz,ds) - v(dz)ds$ is the compensated measure of $N$. This SDE (17) is formally solved as

$$X_t^* = \mathbb{E}\left[X_t^*\right] + Y_t \text{ for } 0 \leq t < 1, \quad (18)$$

where

$$Y_t := \int_0^t \exp\left(-\int_s^t (r+A_\tau)d\tau\right)\int_0^{+\infty} z\tilde{N}(dz,ds) \quad (19)$$

and the mean of $X_t^*$ given by

$$\mathbb{E}\left[X_t^*\right] = x_0 \exp\left(-\int_0^t (r+A_s)ds\right) + \int_0^t \exp\left(-\int_s^t (r+A_\tau)d\tau\right)(-B_s + M_1)ds. \quad (20)$$

The process $Y_t$ for $0 \leq t < 1$ is a local martingale. In the next section, we show that $X_t^*$ satisfies the terminal condition a.s.

**Remark 1** The ODE of $\mathbb{E}\left[X_t^*\right]$ for $0 < t < 1$ is given by

$$\frac{d}{dt}\mathbb{E}\left[X_t^*\right] = -(r+A_t)\mathbb{E}\left[X_t^*\right] - B_t + M_1. \quad (21)$$

The governing linear ODEs of higher-order moments can be obtained by considering the time evolution of suitable monomial of $X^*$. For example, we have for $0 < t < 1$

$$\frac{d}{dt}\mathbb{E}\left[\left(X_t^*\right)^2\right] = -2(r+A_t)\mathbb{E}\left[\left(X_t^*\right)^2\right] + 2(-B_t + M_1)\mathbb{E}\left[X_t^*\right] + M_2. \quad (22)$$

### III. ANALYSIS OF THE STOCHASTIC BRIDGE

This section shows that the obtained process (18) satisfies the terminal condition (2) a.s. To this end, we first show the convergence of the mean (20) and then use it to show the a.s. convergence of (18) as $t \to 1-$.

#### A. Convergence of the expectation

We prove **Proposition 1** which is important on its own and is essential to obtain the pathwise convergence (**Proposition 2**) in the next sub-section. The proof is technical and is long; a more detailed version of the proof is given in the appendix.

**Proposition 1** It follows that

$$\lim_{t \to 1-} \mathbb{E}\left[X_t^*\right] = \hat{x}. \quad (23)$$

**Proof**

By (15), an elementary calculation yields

$$\exp\left(-\int_s^t (r+A_\tau)d\tau\right) = \exp(-r(t-s))\frac{1-\exp(-2r(1-t))}{1-\exp(-2r(1-s))}. \quad (24)$$

For the sake of brevity, we set

$$K_t := \exp(-rt-r)\left(1-\exp(-2r(1-t))\right). \quad (25)$$

We use a series of technical results presented below. Firstly,

$$\int_0^t \frac{2r\hat{x}\exp\left(-\int_s^t (r+A_\tau)d\tau\right)}{\exp(r(1-s))\left(1-\exp(-2r(1-s))\right)}ds = 2r\hat{x}K_t I_{1,t}, \quad (26)$$

with

$$I_{1,t} := \int_0^t \frac{\exp(2rs)}{\left(1-\exp(-2r(1-s))\right)^2}ds$$

$$= \frac{\exp(4r)}{2r}\left[\frac{1}{\exp(2r)-\exp(2rt)} - \frac{1}{\exp(2r)-1}\right]. \quad (27)$$

Secondly,

$$\int_0^t \frac{2M_1\left(1-\exp(-r(1-s))\right)\exp\left(-\int_s^t (r+A_\tau)d\tau\right)}{\exp(r(1-s))\left(1-\exp(-2r(1-s))\right)}ds \quad (28)$$

$$= 2M_1 K_t I_{2,t}$$

with

$$I_{2,t} := \int_0^t \frac{\exp(2rs)\left(1-\exp(-r(1-s))\right)}{\left(1-\exp(-2r(1-s))\right)^2}ds$$

$$= \frac{\exp(3r)}{r}\left[\frac{1}{4\exp(r)}\ln\left(\frac{\exp(r)+u}{\exp(r)-u}\right) + \frac{1}{2}\frac{1}{\exp(r)+u}\right]_1^{\exp(rt)}. \quad (29)$$

Finally,

$$\int_0^t \exp\left(-\int_s^t (r+A_\tau)\,d\tau\right) M_1 ds = M_1 K_t \exp(r) I_{3,t} \quad (30)$$

with

$$I_{3,t} \equiv \int_0^t \frac{\exp(rs)}{1-\exp(-2r(1-s))} ds$$
$$= \frac{\exp(r)}{2r} \left[\ln\left(\frac{\exp(r)+u}{\exp(r)-u}\right)\right]_1^{\exp(rt)}. \quad (31)$$

We then obtain the representation of the mean as

$$\mathbb{E}[X_t^*] = x_0 \exp(-rt) \frac{1-\exp(-2r(1-t))}{1-\exp(-rt)}$$
$$+ K_t \left(2r\hat{x} I_{1,t} - 2M_1 I_{2,t} + M_1 \exp(r) I_{3,t}\right) \quad (32)$$

for $0 \leq t < 1$. By taking the limit $t \to 1-$ using (27), (29), and (31), we obtain the following limits: $\lim_{t\to 1-} 2r\hat{x} K_t I_{1,t} = \hat{x}$ and $\lim_{t\to 1-} K_t I_{2,t} = \lim_{t\to 1-} K_t I_{3,t} = 0$. Consequently, (23) holds.

□

### B. Pathwise Convergence

We now show that the process $X^*$ satisfies the terminal constraint of the OU bridge (18).

**Proposition 2** *It follows that.*

$$\lim_{t\to 1-} X_t^* = \hat{x} \quad a.s. \quad (33)$$

*Proof*

By **Proposition 1** and (18), it suffices to prove $\mathbb{E}[Y_t^2] \to 0$ as $t \to 1-$. Indeed, this together with (19) written in terms of the compensated Poisson measure shows $(Y_t)_{0 \leq t < 1}$ is an $L^2$-bounded martingale and thus martingale convergence theorem shows $Y_t \to 0$ a.s. (e.g., Theorem 27.3 of [21]). For $t$ close to the terminal time 1, $\exp(-rt)\{1 - \exp(-2r(1-t))\} = O(1-t)$. By the isometry (e.g., Problem 9.4 of [22]), we have for $0 \leq t < 1$,

$$\mathbb{E}\left[\left(\int_0^t \frac{\exp(rs)}{1-\exp(-2r(1-s))} \int_0^{+\infty} z\tilde{N}(dz,ds)\right)^2\right]$$
$$= M_2 \mathbb{E}\left[\int_0^t \left(\frac{\exp(rs)}{1-\exp(-2r(1-s))}\right)^2 ds\right] \quad , (34)$$
$$\leq M_2 \exp(2r) \int_0^t \left(\frac{1}{1-\exp(-2r(1-s))}\right)^2 ds$$

which is $O\left((1-t)^{-1}\right)$ as $t \to 1-$. Hence, we obtain, by (34),

$$\mathbb{E}[Y_t^2] = O\left((1-t)^2\right) \times O\left((1-t)^{-1}\right) \to 0 \quad (35)$$

as $t \to 1-$. Therefore, we obtain $Y_t \to 0$ a.s. and together with **Proposition 1**, we arrive at (23).

□

*Remark 2* The theoretical approach in this section carries over to jump processes having bounded variations, but it is not straightforward to extend it to the case driven by jumps of unbounded variation. In particular, the first moment of jumps is not well-defined for the case of unbounded variation. Importantly, the positivity of jumps is essential for the SE bridge presented later.

### C. Self-exciting Bridge

We now extend the results to study, what we call the SE bridge. We consider an SDE with the initial condition $x_0 \in \mathbb{R}$ and for $t > 0$, given by

$$dX_t = -rX_t dt + \int_0^{+\infty} \int_0^{\max\{1+pX_{t-},0\}} zM(dz,dw,dt) \quad (36)$$

with $p \geq 0$ and a Poisson random measure $M$ whose compensator is $v(dz)dwdt$ for the same Lévy measure $v$ on $(0,+\infty)$ as in that for the OU process. Its self-exciting feature is captured by the frequency of positive jumps proportional to the current size of the process. The SDE (36) admits a unique strong càdlàg solution by Theorems 2.3 (i), 3.1, and 2.8 (i) of [23]. The SDE (36) reduces to (1) if $p = 0$. The process $X$ governed by (36) remains non-negative if $X_0 \geq 0$ and hence the max operator in (36) is superficious, but we leave it for the SE bridge constructed below, to make sense.

Now we consider the minimization problem (4) as in the OU case where the controlled process is given by

$$dX_t^{(u)} = \left(-rX_t^{(u)} + u_t\right) dt + \int_0^{+\infty} \int_0^{\max\{1+pX_{t-}^{(u)},0\}} zM(dz,dw,dt) \quad (37)$$

for $t > 0$ and $X_0^{(u)} = x_0$. The objective is again to compute the optimal control to minimize (4) for each $F > 0$.

Fix $F > 0$. The solution to the LQ control becomes $\Phi(0,x_0)$ where $\Phi$ solves the HJB equation

$$\frac{\partial \Phi}{\partial t} - rx \frac{\partial \Phi}{\partial x} + \inf_{u \in \mathbb{R}} \left\{ u \frac{\partial \Phi}{\partial x} + \frac{1}{2} u^2 \right\}$$
$$+ \max\{1+px,0\} \int_0^{+\infty} (\Phi(t,x+z) - \Phi(t,x)) v(dz) = 0 \quad (38)$$

for $0 \leq t < 1$ and $x \in \mathbb{R}$ subject to the terminal condition (6). The solution to (38) fails to be of the quadratic form (7), but a modification of (38) with $\max\{1+px,0\}$ replaced by $1+px$ can be solved analytically. The solution of the modified HJB equation is of the form (7), in terms of the Riccati equation of the form (8)-(11) with $r$ replaced by $R$ and $M_1$ in (9) by $m$, where $R := r - pM_1$ and $m := M_1 + pM_2/2$. This modified version is expected to approximate the exact solution of (38) when $p \geq 0$ is small. After taking the limit $F \to +\infty$, we arrive at the coefficients $A_t, B_t$ with $r$ replaced by $R$ and $M_1$ in (9) by $m$, which enables us to obtain the LQ regulator

of the form (12) explicitly. We then obtain (the approximation of the) SE bridge given by

$$dX_t^* = \left(-(R+A_t)X_t^* - B_t\right)dt$$
$$+ \int_0^{+\infty} \int_0^{\max\{1+pX_{t-}^*, 0\}} zM(dz, dw, dt). \quad (39)$$

**Remark 3** If the controlled process $X^*$ is larger than $-p^{-1}$ throughout time interval $[0,1]$ a.s., then the max operator in (36) becomes superficious and the constructive argument of the bridge applies without approximations. This condition is non-trivial *a priori* due to the possible positivity of $B_t$ near the terminal time. We computationally see in the next section that this pathwise non-negative property is satisfied when $p > 0$ is small.

## IV. COMPUTATIONAL EXPERIMENTS

### A. Numerical Method

We verify the analytical results obtained in the previous sections via a Monte-Carlo simulation using the classical Euler–Maruyama discretization scheme with the fixed step size $\Delta t = 1/200,000$ and 200,000 sample paths. Random numbers for the jump noises are generated by the Mersenne Twister [24]. The ODEs governing the first- and second-order moments ((21) and (22) in **Remark 1**) of the bridges are computed numerically by the classical forward finite-difference Euler discretization scheme with the same step size, sufficiently small to achieve accurate approximation. Below, we use two cases for the Lévy measure $v$: finite activity case (a compound Poisson case) and infinite activity case.

### B. Computational Results

We set $x_0 = \hat{x} = 0$. **Fig. 1** shows sample paths of the OU ($p=0$) and SE bridges ($p=2$) driven by a compound Poisson process with $v(dz) = 2\exp(-50z)dz$ for $z > 0$ and the reversion rate $r = 10$. As observed in **Fig. 1**, the SE bridge involves more clustered jumps due to its self-exciting mechanism. We evaluated the error of the terminal condition $X_1^* = 0$ using $\text{Err} := \sum_{n=1}^{200,000} |X_{1,n}^*|/200,000$ where $X_{1,n}^*$ is the terminal value of the $n$ th sample path. We obtained Err = $1.91 \times 10^{-5}$ for the OU bridge and is Err = $2.01 \times 10^{-5}$ for SE bridge, supporting the a.s. convergence as verified in **Proposition 2**. Moreover, with $p = 2$, all 200,000 sample paths of the SE bridge are valued in $(-p^{-1}, +\infty)$, supporting the conjecture in **Remark 2**. However, our experiment also showed that this condition fails when we select a larger value of $p$ such as $p = 6$. **Figs. 2-3** plot the mean and variance of the OU and SE bridges computed using the finite-difference and Monte-Carlo methods. The differences are almost unrecognizable, confirming the obtained analytical results. The higher mean and variance in the SE bridge compared to the OU counterpart is due to its self-exciting behavior.

We also apply the SE bridge driven by a subordinator with infinite activities for the problem of streamflow regulation. We assume a one-month time horizon with $T = 1$ (month) or equivalently 30 (days). The parameter values we use here were calibrated using the discharge time series data at the Kamiyasuda Station of the Gono River, Hiroshima Prefecture, Japan. This station is located at the upstream of the Haidzuka Dam for flood control and water supply. The dam is also serving as a habitat for the landlocked Ayu Sweetfish *Plecoglossus altivelis altivelis* as a major fishery resource of Gono River [25]. Thus, controlling the inflow discharge to the dam is vital for regional resource, environmental, and disaster management. The one-year hourly discharge data at the station (**Fig. 4**) was obtained from Water Information System by Ministry of Land, Infrastructure, Transport and Tourism, Japan (http://www1.river.go.jp/). Based on the time series data, $r$ is calibrated by least-squares regression of the autocorrelation function to be 15.8 (1/month). By the moment matching method under a stationary condition [4] fitted to a tempered-stable type Lévy measure, we obtained $v(\mathrm{d}z) = 3.23\exp(-0.031z)z^{-1.87}\mathrm{d}z$ and $p = 0.14$.

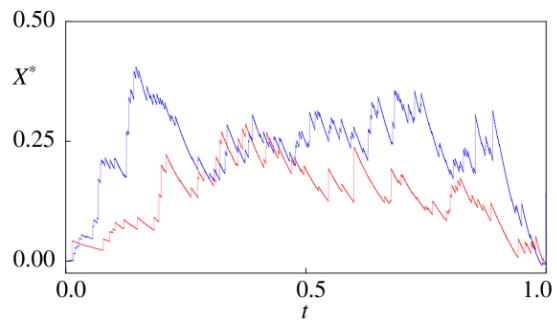

Figure 1. Sample paths of the OU (red) and SE (blue) bridges driven by a compound Poisson process.

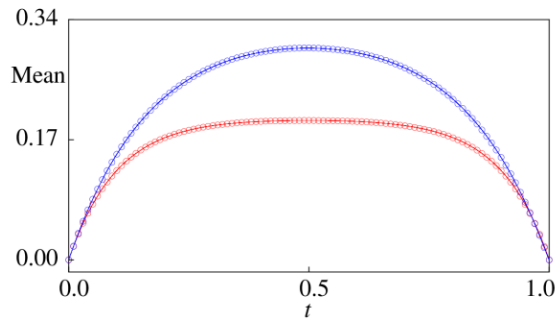

Figure 2. Comparison of the means of the OU (red) and SE bridges (blue) between theoretical finite-difference (lines) and Monte-Carlo results (circles).

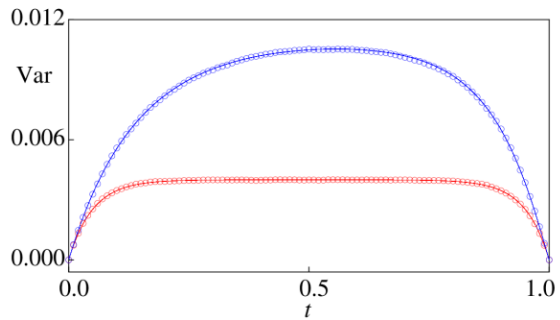

Figure 3. Comparison of the variances (Var) of the OU (red) and SE bridges (blue) between theoretical finite difference (lines) and Monte-Carlo results (circles).

We consider a scenario where the discharge needs to be maintained to be moderately high at the terminal time ($\hat{x} = 4$ or $8$ (m³/s)) by implementing a controllable weir at the observation station. For Monte-Carlo simulation, the infinite-activity jump noise is generated using the acceptance-rejection method (Algorithm 0 of [26]). **Fig. 5** plots computed sample paths and the comparison of the mean between the finite-difference and Monte-Carlo results, justifying the obtained theory even for an infinite activity case. All 200,000 sample paths are valued in $(0,+\infty)$, and hence, in view of **Remark 2**, the approximation of (38) is practical with this choice of $p$. Err of the terminal condition is almost negligible with $5.11\times10^{-4}$ (m³/s) for $\hat{x} = 4$ (m³/s) and is $5.65\times10^{-4}$ (m³/s) for $\hat{x} = 8$ (m³/s).

**Fig. 6** shows the coefficients $A, B$ of the OU bridge in **Fig. 2**. Similarly, **Fig. 7** plots $A, B$ of the SE bridge in **Fig. 5** for $\hat{x} = 8$ (m³/s). As shown in these figures, both A and B stay near zero in the early state and grow sharply as it approaches the terminal time. In particular, from **Fig. 7**, the implication of this is that the streamflow should be controlled intensively near the terminal time so that the designed terminal condition is met with the energetic optimality. This kind of control strategies would be important especially for the flood control where the rapid decision-making is required. As a function of $t$, the mean first increases, flattens and then increases to satisfy the terminal condition. In other words, the obtained SE bridge has a turnpike property [e.g., 28]. The streamflow would be stabilized in most of the time horizon.

In reality, the capacity of control is limited and it makes sense to consider a problem with regularization or saturation. However, this extension is out of scope of this paper, because our objective in this paper is to obtain an explicit solution which is easy to implement. Importantly, it is of mathematical value to obtain a bridge satisfying the terminal condition. It is also noted that this paper also studies the case with relaxed condition with finite penalizing parameter $F$. Nonetheless, the extension with regularization or saturation is indeed important and hence we leave this extension a future work. This requires a completely different, numerical approach because the solution to the corresponding HJB equation no longer admits a quadratic form.

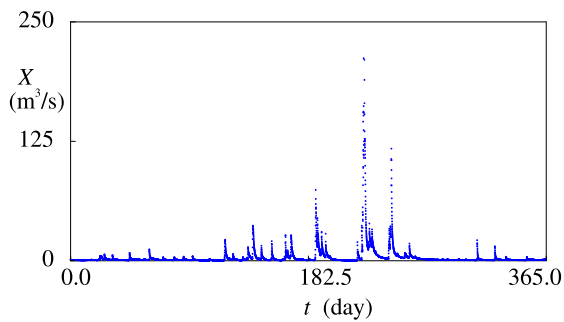

Figure 4. Hourly discharge time series data at Kamiyasuda Station.

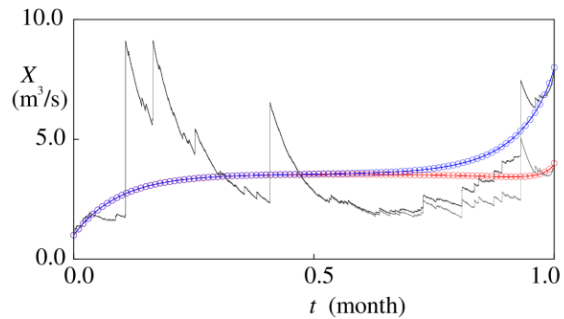

Figure 5. Sample paths (black: $\hat{x} = 8$, grey: $\hat{x} = 4$) and means of the SE bridge between finite difference (blue line: $\hat{x} = 8$, red line: $\hat{x} = 4$) and Monte-Carlo results (blue circles: $\hat{x} = 8$, red circles: $\hat{x} = 4$).

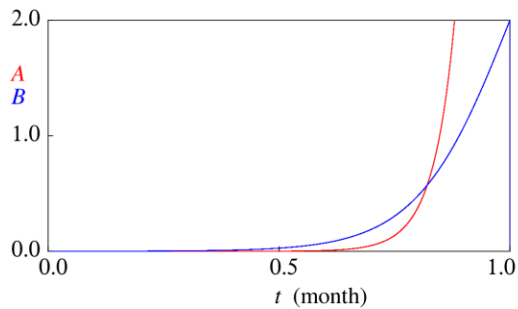

Figure 6. The coefficients $A, B$ for the OU bridge in Figure 2.

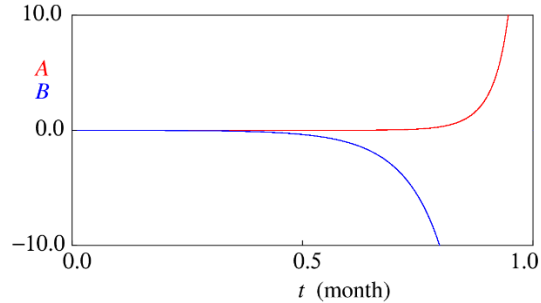

Figure 7. The coefficients $A, B$ for the SE bridge for $\hat{x} = 8$ in Figure 5.

## V. CONCLUSION

This paper studied subordinator-driven stochastic bridges, written in terms of a certain limit of the solution to LQ control problems. The terminal condition was confirmed to be satisfied both analytically and numerically. There are several venues for future work. First, it is of interest to consider additional constraints on the maximum and integrated value of the path until the terminal time and develop bridges satisfying these conditions. This extension is applicable to, for example, the control of a dam-reservoir system where the water storage level needs to be controlled while satisfying several operational constraints. Another direction is to consider stochastic bridges driven by more general processes. The mathematical framework to construct stochastic bridges in this paper is extendible to, for example, multi-dimensional and long-memory cases, but the resulting Riccati equation will not admit analytical solutions. For these extensions, it is vital to develop efficient computational methods.

**Appendix of "A Jump Ornstein–Uhlenbeck Bridge Based on Energy-optimal Control and Its Self-exciting Extension" by H. Yoshioka and K. Yamazaki**

In this appendix, we prove **Proposition 1** in the main text with more technical details.

**Proposition 1** *It follows that*

$$\lim_{t \to 1-} \mathbb{E}\left[X_t^*\right] = \hat{x}. \qquad (40)$$

**Proof**

First of all, we have

$$\mathbb{E}\left[X_t^*\right] = x_0 \exp\left(-\int_0^t (r + A_s) \mathrm{d}s\right) + \int_0^t \exp\left(-\int_s^t (r + A_\tau) \mathrm{d}\tau\right)(-B_s + M_1) \mathrm{d}s \text{ for } 0 \le t < 1. \qquad (41)$$

We prove (40) using this representation of $\mathbb{E}\left[X_t^*\right]$. By an elementary calculation, we obtain

$$\begin{aligned}
&\exp\left(-\int_s^t (r + A_\tau) \mathrm{d}\tau\right) \\
&= \exp(-r(t-s))\exp\left(-\int_s^t A_\tau \mathrm{d}\tau\right) \\
&= \exp(-r(t-s))\exp\left(-\int_s^t \frac{2r}{\exp(2r(1-\tau))-1} \mathrm{d}\tau\right) \\
&= \exp(-r(t-s))\exp\left(-\int_s^t \frac{2r\exp(2r\tau)}{\exp(2r)-\exp(2r\tau)} \mathrm{d}\tau\right) \\
&= \exp(-r(t-s))\exp\left(\left[\ln(\exp(r)-\exp(2r\tau))\right]_s^t\right) \\
&= \exp(-r(t-s))\exp\left(\ln\left(\frac{\exp(2r)-\exp(2rt)}{\exp(2r)-\exp(2rs)}\right)\right) \\
&= \exp(-r(t-s))\exp\left(\ln\left(\frac{1-\exp(-2r(1-t))}{1-\exp(-2r(1-s))}\right)\right) \\
&= \exp(-r(t-s))\frac{1-\exp(-2r(1-t))}{1-\exp(-2r(1-s))}
\end{aligned} \qquad (42)$$

For brevity, we set

$$K_t := \exp(-rt-r)(1-\exp(-2r(1-t))). \qquad (43)$$

Then, we have

$$\begin{aligned}
-B_s &= \frac{2r}{\exp(r(1-s))(1-\exp(-2r(1-s)))}\left(\hat{X} - \frac{M}{r}(1-\exp(-r(1-s)))\right) \\
&= \frac{2r\hat{X}}{\exp(r(1-s))(1-\exp(-2r(1-s)))} - \frac{2M(1-\exp(-r(1-s)))}{\exp(r(1-s))(1-\exp(-2r(1-s)))}
\end{aligned}. \qquad (44)$$

We use a series of technical results presented below. Firstly,

$$\int_0^t \exp\left(-\int_s^t (r+A_\tau)\,d\tau\right) \frac{2r\hat{x}}{\exp(r(1-s))(1-\exp(-2r(1-s)))}\,ds$$

$$= 2r\hat{x}\int_0^t \exp(-r(t-s)) \frac{1-\exp(-2r(1-t))}{1-\exp(-2r(1-s))} \frac{1}{\exp(r(1-s))(1-\exp(-2r(1-s)))}\,ds$$

$$= 2r\hat{x}\exp(-rt)(1-\exp(-2r(1-t)))\int_0^t \exp(rs) \frac{1}{\exp(r(1-s))(1-\exp(-2r(1-s)))^2}\,ds \quad (45)$$

$$= 2r\hat{x}K_t \int_0^t \frac{\exp(2rs)}{(1-\exp(-2r(1-s)))^2}\,ds$$

$$= 2r\hat{x}K_t I_{1,t}$$

with

$$I_{1,t} := \int_0^t \frac{\exp(2rs)}{(1-\exp(-2r(1-s)))^2}\,ds$$

$$= \frac{1}{2r}\int_0^t \frac{1}{(1-\exp(-2r(1-s)))^2}\,d(\exp(2rs))$$

$$= \frac{1}{2r}\int_1^{\exp(2rt)} \frac{1}{(1-u\exp(-2r))^2}\,du \qquad (46)$$

$$= \frac{\exp(4r)}{2r}\int_1^{\exp(2rt)} \frac{1}{(\exp(2r)-u)^2}\,du$$

$$= \frac{\exp(4r)}{2r}\left[\frac{1}{\exp(2r)-u}\right]_1^{\exp(2rt)}$$

$$= \frac{\exp(4r)}{2r}\left[\frac{1}{\exp(2r)-\exp(2rt)} - \frac{1}{\exp(2r)-1}\right]$$

Secondly, we have

$$\int_0^t \exp\left(-\int_s^t (r+A_\tau)\,d\tau\right)\left(-\frac{2M_1(1-\exp(-r(1-s)))}{\exp(r(1-s))(1-\exp(-2r(1-s)))}\right)ds$$

$$= -2M_1 \int_0^t \exp(-r(t-s)) \frac{1-\exp(-2r(1-t))}{1-\exp(-2r(1-s))} \frac{1-\exp(-r(1-s))}{\exp(r(1-s))(1-\exp(-2r(1-s)))}\,ds \quad (47)$$

$$= -2M_1 K_t \int_0^t \frac{\exp(2rs)(1-\exp(-r(1-s)))}{(1-\exp(-2r(1-s)))^2}\,ds$$

$$= -2M_1 K_t I_{2,t}$$

with

$$I_{2,t} := \int_0^t \frac{\exp(2rs)\left(1-\exp(-r(1-s))\right)}{\left(1-\exp(-2r(1-s))\right)^2} ds$$

$$= \frac{1}{r}\int_0^t \frac{\exp(rs)\left(1-\exp(-r(1-s))\right)}{\left(1-\exp(-2r(1-s))\right)^2} d(\exp(rs))$$

$$= \frac{1}{r}\int_1^{\exp(rt)} \frac{u\left(1-u\exp(-r)\right)}{\left(1-u^2\exp(-2r)\right)^2} du$$

$$= \frac{\exp(4r)}{r}\int_1^{\exp(rt)} \frac{u\left(1-u\exp(-r)\right)}{\left(\exp(2r)-u^2\right)^2} du$$

$$= \frac{\exp(3r)}{r}\int_1^{\exp(rt)} \frac{u\left(\exp(r)-u\right)}{\left(\exp(2r)-u^2\right)^2} du$$

$$= \frac{\exp(3r)}{r}\int_1^{\exp(rt)} \frac{u\left(\exp(r)-u\right)}{\left(\exp(r)-u\right)^2\left(\exp(r)+u\right)^2} du$$

$$= \frac{\exp(3r)}{r}\int_1^{\exp(rt)} \frac{u}{\left(\exp(r)-u\right)\left(\exp(r)+u\right)^2} du$$

$$= \frac{\exp(3r)}{r}\int_1^{\exp(rt)} \frac{u-\exp(r)+\exp(r)}{\left(\exp(r)-u\right)\left(\exp(r)+u\right)^2} du$$

$$= \frac{\exp(3r)}{r}\int_1^{\exp(rt)} \left\{\frac{\exp(r)}{\left(\exp(r)-u\right)\left(\exp(r)+u\right)^2} - \frac{1}{\left(\exp(r)+u\right)^2}\right\} du \qquad (48)$$

and

$$\frac{\exp(r)}{\left(\exp(r)-u\right)\left(\exp(r)+u\right)^2} - \frac{1}{\left(\exp(r)+u\right)^2}$$

$$= \frac{1}{2}\frac{1}{\exp(r)+u}\left(\frac{1}{\exp(r)-u}+\frac{1}{\exp(r)+u}\right) - \frac{1}{\left(\exp(r)+u\right)^2}, \qquad (49)$$

$$= \frac{1}{2}\frac{1}{\exp(r)+u}\frac{1}{\exp(r)-u} - \frac{1}{2}\frac{1}{\left(\exp(r)+u\right)^2}$$

$$= \frac{1}{4\exp(r)}\left(\frac{1}{\exp(r)+u}+\frac{1}{\exp(r)-u}\right) - \frac{1}{2}\frac{1}{\left(\exp(r)+u\right)^2}$$

leading to

$$I_{2,t} = \frac{\exp(3r)}{r}\int_1^{\exp(rt)} \left\{\frac{\exp(r)}{\left(\exp(r)-u\right)\left(\exp(r)+u\right)^2} - \frac{1}{\left(\exp(r)+u\right)^2}\right\} du$$

$$= \frac{\exp(3r)}{r}\int_1^{\exp(rt)} \left\{\frac{1}{4\exp(r)}\left(\frac{1}{\exp(r)+u}+\frac{1}{\exp(r)-u}\right) - \frac{1}{2}\frac{1}{\left(\exp(r)+u\right)^2}\right\} du . \qquad (50)$$

$$= \frac{\exp(3r)}{r}\left[\frac{1}{4\exp(r)}\ln\left(\frac{\exp(r)+u}{\exp(r)-u}\right) + \frac{1}{2}\frac{1}{\exp(r)+u}\right]_1^{\exp(rt)}$$

Finally, we have

$$\int_0^t \exp\left(-\int_s^t (r+A_\tau)\,d\tau\right) M_1\,ds = M_1 \int_0^t \exp(-r(t-s)) \frac{1-\exp(-2r(1-t))}{1-\exp(-2r(1-s))}\,ds$$
$$= M_1 \exp(-rt)(1-\exp(-2r(1-t))) I_{3,t} \qquad (51)$$
$$= M_1 K_t \exp(r) I_{3,t}$$

with

$$\begin{aligned}
I_{3,t} &:= \int_0^t \frac{\exp(rs)}{1-\exp(-2r(1-s))}\,ds \\
&= \frac{1}{r}\int_0^t \frac{1}{1-\exp(-2r(1-s))}\,d(\exp(rs)) \\
&= \frac{1}{r}\int_0^t \frac{\exp(2r)}{\exp(2r)-\exp(2rs)}\,d(\exp(rs)) \\
&= \frac{\exp(2r)}{r}\int_1^{\exp(rt)} \frac{1}{\exp(2r)-u^2}\,du \\
&= \frac{\exp(2r)}{r}\frac{1}{2\exp(r)}\int_1^{\exp(rt)} \left(\frac{1}{\exp(r)+u}+\frac{1}{\exp(r)-u}\right)du \\
&= \frac{\exp(r)}{2r}\left[\ln\left(\frac{\exp(r)+u}{\exp(r)-u}\right)\right]_1^{\exp(rt)}
\end{aligned} \qquad (52)$$

Consequently, we obtain

$$\mathbb{E}[X_t] = x_0 \exp(-rt)\frac{1-\exp(-2r(1-t))}{1-\exp(-2r)} \qquad \text{for } 0<t<1. \qquad (53)$$
$$+ K_t\left(2r\hat{x} I_{1,t} - 2M I_{2,t} + M_1 \exp(r) I_{3,t}\right)$$

By taking the limit $t \to 1-$ using (46), (50), (52), we obtain the following limits. Firstly, we obtain

$$\lim_{t\to 1-} 2r\hat{x} K_t I_{1,t} = \hat{x}. \qquad (54)$$

Indeed, we have

$$\begin{aligned}
&\lim_{t\to 1-} 2r\hat{x} K_t I_{1,t} \\
&= \lim_{t\to 1-} 2r\hat{x} K_t \frac{\exp(4r)}{2r}\left[\frac{1}{\exp(2r)-\exp(2rt)} - \frac{1}{\exp(2r)-1}\right] \\
&= \hat{x} \lim_{t\to 1-} \exp(-rt-r)(1-\exp(-2r(1-t)))\exp(4r)\left[\frac{1}{\exp(2r)-\exp(2rt)} - \frac{1}{\exp(2r)-1}\right] \\
&= \hat{x} \lim_{t\to 1-} \exp(3r-rt)\frac{1-\exp(-2r(1-t))}{\exp(2r)-\exp(2rt)} \\
&= \hat{x} \lim_{t\to 1-} \exp(r-rt)\frac{\exp(2r)-\exp(2rt)}{\exp(2r)-\exp(2rt)} \\
&= \hat{x} \lim_{t\to 1-} \exp(r-rt) \\
&= \hat{x}
\end{aligned} \qquad (55)$$

Similarity, we obtain

$$\lim_{t \to 1-} K_t I_{2,t} = \lim_{t \to 1-} K_t I_{3,t} = 0. \tag{56}$$

Indeed, we have

$$\lim_{t \to 1-} K_t I_{2,t}$$

$$= \lim_{t \to 1-} \exp(-rt - r)(1 - \exp(-2r(1-t))) \frac{\exp(3r)}{r} \left[ \frac{1}{4\exp(r)} \ln\left(\frac{\exp(r) + u}{\exp(r) - u}\right) + \frac{1}{2} \frac{1}{\exp(r) + u} \right]_1^{\exp(rt)}$$

$$= \lim_{t \to 1-} \exp(-rt)(1 - \exp(-2r(1-t))) \frac{\exp(2r)}{r} \left[ \frac{1}{4\exp(r)} \ln\left(\frac{\exp(r) + u}{\exp(r) - u}\right) \right]_1^{\exp(rt)}$$

$$= \frac{\exp(r)}{4r} \lim_{t \to 1-} \exp(-rt)(1 - \exp(-2r(1-t))) \left[ \ln\left(\frac{\exp(r) + \exp(rt)}{\exp(r) - \exp(rt)}\right) - \ln\left(\frac{\exp(r) + 1}{\exp(r) - 1}\right) \right] \tag{57}$$

$$= \frac{\exp(r)}{4r} \lim_{t \to 1-} \exp(-rt)(1 - \exp(-2r(1-t))) \left[ \ln\left(\frac{\exp(r) + \exp(rt)}{\exp(r) + 1}\right) + \ln\left(\frac{\exp(r) - 1}{\exp(r) - \exp(rt)}\right) \right]$$

$$= \frac{\exp(r)}{4r} \lim_{t \to 1-} \exp(-rt)(1 - \exp(-2r(1-t))) \ln\left(\frac{\exp(r) - 1}{\exp(r) - \exp(rt)}\right)$$

$$= \lim_{t \to 1-} \left\{ O(1-t) \times O\left(\ln\left(\frac{1}{1-t}\right)\right) \right\}$$

$$= 0$$

and

$$\lim_{t \to 1-} K_t I_{3,t} = \lim_{t \to 1-} \exp(-rt - r)(1 - \exp(-2r(1-t))) \frac{\exp(r)}{2r} \left[ \ln\left(\frac{\exp(r) + u}{\exp(r) - u}\right) \right]_1^{\exp(rt)}$$

$$= \lim_{t \to 1-} \left\{ O(1-t) \times O\left(\ln\left(\frac{1}{1-t}\right)\right) \right\} \tag{58}$$

$$= 0$$

Consequently, we obtain the desired result (40).

□